\documentclass[a4paper,10pt]{article}
\usepackage[utf8]{inputenc}
\usepackage{amsmath}
\usepackage{amssymb}

\usepackage{tikz-cd}
\usepackage{showlabels}
\input V3,95
\newtheorem{thm}{Theorem}[section]
 \newtheorem{cor}[thm]{Corollary}
 \newtheorem{lem}[thm]{Lemma}
 \newtheorem{prop}[thm]{Proposition}

 \newtheorem{rem}[thm]{Remark}

 \newtheorem{defn}[thm]{Definition}
 \newtheorem{ex}[thm]{Example}
 \numberwithin{equation}{section}

\newcommand{\Der}{{\rm Der}}

\newcommand{\Hom}{{\rm Hom}}

\newcommand{\ext}{{\rm Ext}}
\newcommand{\im}{{\rm im}}

\title{Finite generation of Andr\'e-Quillen (co-)homology of F-finite algebras}
\author{Cristodor Ionescu\\Simion Stoilow Institute of Mathematics \\ of the Romanian Academy\\
  P.O. Box 1-764\\ RO 014700 Bucharest\\  Romania
  \\
email: cristodor.ionescu@imar.ro}

\begin{document}

\maketitle

\begin{abstract}
We prove that the Andr\'e-Quillen homology and cohomology modules of F-finite $\mathbb{Z}_{(p)}$-algebras are finitely generated.
\end{abstract}

\maketitle
\section{Introduction and preliminaries}
Andr\'e-Quillen homology was invented independently by M. Andr\'e and D. Quillen (see \cite{A1}, \cite{Qu} and  \cite{A2}). The theory proved its importance over the years, for example by characterizing several classes of noetherian rings and morphisms between such rings (see \cite{A2}, \cite{Qu} and \cite{Iy}). The finite generation of certain Andr\'e-Quillen homology modules appears in many situations, maybe the most well-known being the proof of the celebrated Andr\'e's Theorem on the localization of formal smoothness (see \cite{A3}). In a recent paper, Dundas and Morrow \cite{DM} proved that for a F-finite algebra $A$, the modules $H_i(C,A,E), i\geq 1$ are finitely generated, where $C$ is any of the rings $\mathbb{Z}, \mathbb{Z}_{(p)}$ or $\mathbb{Z}/p^e\mathbb{Z}$ and $E$ is a finitely generated $A$-module. Actually they proved more, namely that also the higher Andr\'e-Quillen homology modules are finitely generated, as their final goal was to prove the finite generation of the Hochschild homology modules. Our purpose is to give another proof of the finite generation of $H_i(C,A,E)$, where $C$ is as above, as well as some similar results concerning the Andr\'e-Quillen cohomology modules.

\par All the rings are commutative, with unit and noetherian. All over the paper  $p>0$ will be a fixed  prime number. We will denote by $\mathbb{F}_p$ the prime field of characteristic $p$ and by $\mathbb{Z}_{(p)}$ the ring of fractions of $\mathbb{Z}$ with denominators not divisible by $p.$ For a ring $A$ containing the field $\mathbb{F}_p,$ the Frobenius morphism of $A$ will mean the ring morphism $F:A\to A, F(a)=a^p.$  

\begin{rem}\label{alg} Let $A$ be a $\mathbb{Z}_{(p)}$-algebra. We have a canonical commutative diagram of ring morphisms
 
\begin{diagram}[size=2.8em]
 \mathbb{Z}   &  \rTo  &   \mathbb{Z}_{(p)} &\rTo & A \\
 \dTo         &       &   \dTo  &   &   \dTo \\
 \mathbb{Z}/p^e\mathbb{Z} &\rTo   &\mathbb{Z}/p\mathbb{Z}\simeq\mathbb{Z}_{(p)}/p\mathbb{Z}_{(p)}\simeq\mathbb{F}_p  &  \rTo  &  A/pA.
\end{diagram}


 There are two possibilities:
\par 1) $pA=0.$ Then $A=A/pA$  contains $\mathbb{F}_p.$
\par 2) $pA\neq 0.$ Then $A$ does not contain $\mathbb{F}_p,$ but $A/pA$ does.\\
\end{rem}
\par Thus in both cases we can say that $A/pA$ contains $\mathbb{F}_p$ and the following definition makes sense and generalizes the usual one (see \cite{LZ} and \cite{DM}):
\begin{defn}\label{zpf}
 A noetherian $\mathbb{Z}_{(p)}$-algebra A is called F-finite, if the Frobenius morphism of $A/pA$ is a finite morphism.
\end{defn}

We will heavily rely on the following result of Gabber:

\begin{thm}\label{gab}{\rm \cite[Rem. 13.6]{Gab}}
 Let $k$ be a field of characteristic $p>0$ and A be a noetherian F-finite k-algebra. Then A is the quotient of a regular F-finite k-algebra.
\end{thm}

\begin{rem}\label{obsalg}
Theorem \ref{gab} does not extend to the case 2) of remark \ref{alg}, that is to algebras not containing a field of characteristic $p$.
\end{rem}

\section{Finiteness of Andr\'e-Quillen (co-)homology}
The next lemma is well-known, but we couldn't find a precise reference.
\begin{lem}\label{sirex}
Let A be a noetherian ring and $E\stackrel{u}{\to} F\stackrel{v}{\to} G$ be an exact sequence of A-modules. If E and G are finitely generated, then 
F is finitely generated.
\end{lem}

\par\noindent\textit{Proof:}  We have the exact sequences 
$$0\to \ker(u)\to E\to E/\ker(u)\cong\im(u)\to 0$$
and
$$0\to \ker(v)\to F\to F/\ker(v)\cong \im(v)\to 0.$$
From the first sequence it follows that $\im(u)=\ker(v)$ is finitely generated and since $\im(v)$ is finitely generated, 
from the second one we obtain the conclusion.

\begin{prop}\label{f-finit}
 Let k be a perfect field of characteristic $p>0$ and A be a noetherian F-finite k-algebra. Then $H_i(k,A,E)$ and $H^i(k,A,E)$ are  finitely generated $A$-modules for
 any $i\geq 0$ and for any finitely generated $A$-module $E.$
 \end{prop}
\par\noindent\textit{Proof:} 
 From \ref{gab} it follows that  there exists a regular ring $R$ and a surjective morphism $R\to A.$ Consider the Jacobi-Zariski exact sequence attached to $k\to R \to A.$ For any $i\geq 1$ we get an exact sequence 
 $$H_i(k,R,E)\to H_i(k,A,E)\to H_i(R,A,E).$$
 From \cite[Prop. IV.55]{A2},  it follows that $H_i(R,A,E)$  is finitely generated. As $k$ is perfect and $R$ is regular, we get that $k\to R$ is a regular morphism and from \cite[Thm. 9.5]{Iy},  we have that $H_i(k,R,E)=0, \forall i\geq 1,$ whence $H_i(k,A,E)$ is finitely generated. For $i=0$ we have that 
 $H_0(k,A,E)\cong\Omega_{A/k}=\Omega_{A/k[A^p]}$ which is finitely generated.\\
  As $R$ is F-finite, the module of differentials $\Omega_{R/k}$ is finitely generated. Since  $H_i(k,R,R)=0$ for all $i\geq 1,$  by \cite[Prop. 3.19]{A2} it follows that $H^j(k,R,E)\simeq \ext^j_R(\Omega_{R/k},E), j\geq 0,$ hence by \cite[Th. 7.36]{Ro} we obtain that $H^j(k,R,E)$ is finitely generated. Now the Jacobi-Zariski exact sequence associated to $k\to R\to A$  is 
$$H^i(R,A,E)\to H^i(k,A,E)\to H^i(k,R,E)$$
and \ref{sirex} ends the proof.

\begin{rem}\label{general} Along the same lines one can prove that, if $k$ is a perfect field (possibly of characteristic zero) and $A$ is a $k$-algebra which is a quotient of a regular ring, then $H_i(k,A,E)$ is a finitely generated $A$-module for any $i\geq 1$ and any finitely generated $A$-module $E.$ 
\end{rem}

\begin{prop}\label{2}
 Let k be a  field of characteristic $p>0,$  A be a noetherian F-finite k-algebra and E a finitely generated A-module. Then:
 \par a) $H_i(k,A,E)$ and $H^i(k,A,E)$ are finitely generated $A$-modules for
 any $i\geq 0, i\neq 1$;
 \par b) If moreover $k$ is F-finite, then $H_1(k,A,E)$ and $H^1(k,A,E)$ are finitely generated $A$-modules.
\end{prop}

\par\noindent\textit{Proof:}  a) Let us consider the Jacobi-Zariski exact sequence in homology associated to the morphisms $\mathbb{F}_p\to k\to A,$ 
namely
$$\ldots\to H_i(\mathbb{F}_p,k,E)\to H_i(\mathbb{F}_p,A,E)\to H_i(k,A,E)\to H_{i-1}(\mathbb{F}_p,k,E)\to\ldots$$
But for $i\geq 2,$ by the separability of $k$ over $\mathbb{F}_p$ and \cite[Prop. VII.11 and VII.4]{A2} we have that $H_i(\mathbb{F}_p,A,E)\cong H_i(k,A,E),$ 
hence by \ref{f-finit} we obtain the assertion. For $i=0$ we have $H_0(k,A,E)=\Omega_{A/k}\otimes_AE=\Omega_{A/k[A^p]}\otimes_A E,$ hence it is finitely generated.\\
Let us now consider the Jacobi-Zariski exact sequence in cohomology associated to the same morphisms as above, that is 
$$\ldots\to H^{i-1}(\mathbb{F}_p,k,E)\to H^i(k,A,E)\to H^i(\mathbb{F}_p,A,E)\to H^i(\mathbb{F}_p,k,E)\to\ldots$$
Again, for $i\geq 2,$ by the separability of $k$ over $\mathbb{F}_p$ and \cite[Prop. VII.4]{A2} we get the isomorphism $H^i(\mathbb{F}_p,A,E)\cong H^i(k,A,E),$ 
hence by \ref{f-finit} we obtain the assertion.
\par\noindent b) We have the exact sequences 
$$0\to H_1(\mathbb{F}_p,A,E)\to H_1(k,A,E)\to \Omega_k\otimes_A E$$
and 
$$\Der_{\mathbb{F}_p}(k,E)\to H^1(k,A,E)\to H^1(\mathbb{F}_p,A,E)\to H^1(\mathbb{F}_p,k,E)=(0).$$
But the hypothesis implies that $\Omega_k$ is finitely generated, hence from \ref{f-finit} and \ref{sirex} we get that $H_1(k,A,E)$ is finitely generated. Moreover $\Der_{\mathbb{F}_p}(k,E)=\Hom_k(\Omega_{k/\mathbb{F}_p},E)$ which is finitely generated, since
$\Omega_{k/\mathbb{F}_p}$ is finitely generated. From the second exact sequence, \ref{sirex} and \ref{f-finit} we obtain the desired conclusion.

\begin{ex}\label{cex}  If $[k:k^p]=\infty$ the module $H_1(k,A,E)$ is not necessarily finitely generated. Indeed, let $k$ be a field such that $[k:k^p]=\infty$ and let $\bar{k}$ be the algebraic closure of $k.$ The Jacobi-Zariski sequence associated to $\mathbb{F}_p\to k\to \bar{k}$ is
$$0\to H_1(k,\bar{k},\bar{k})\to\Omega_k\otimes\bar{k}\to\Omega_{\bar{k}}\to\Omega_{\bar{k}/k}\to 0.$$
Then $\Omega_{\bar{k}}$ is finitely generated and if $H_1(k,\bar{k},\bar{k})$ is finitely generated, by \ref{sirex} it follows that $\Omega_k$ is finitely generated, contradicting the hypothesis that $[k:k^p]=\infty.$
 \end{ex}

\begin{cor}\label{z}
 Let p be a prime number, A be a noetherian F-finite $\mathbb{F}_p$-algebra and E a finitely generated A-module. Then for any $i\geq 0$ and  any $e\geq 1,$ the $A$-modules $H_i(\mathbb{Z},A,E)$,  $H_i(\mathbb{Z}/p^e\mathbb{Z},A,E),$ $H^i(\mathbb{Z},A,E)$ and $H^i(\mathbb{Z}/p^e\mathbb{Z},A,E)$ are  finitely generated.
 \end{cor}

\par\noindent\textit{Proof:}  Considering the Jacobi-Zariski exact sequence associated to the morphisms 
$\mathbb{Z}\to\mathbb{Z}/p\mathbb{Z}=\mathbb{F}_p\to A$ we get an exact sequence 
$$H_i(\mathbb{Z},\mathbb{F}_p,E)\to H_i(\mathbb{Z},A,E)\to H_i(\mathbb{F}_p,A,E).$$ By \ref{f-finit} we know that 
$H_i(\mathbb{F}_p,A,E)$ is finitely generated and by \cite{A2}, IV.55 we obtain that $H_i(\mathbb{Z},\mathbb{F}_p,E)$ is 
finitely generated. Now \ref{sirex} tells us that $H_i(\mathbb{Z},A,E)$ is finitely generated. For the second assertion let us remark first that we can assume $e\geq 2.$ Consider the morphisms $\mathbb{Z}/p^e\mathbb{Z}\to\mathbb{Z}/p\mathbb{Z}=\mathbb{F}_p\to A$ and the Jacobi-Zariski associated sequence 
$$H_i(\mathbb{Z}/p^e\mathbb{Z},\mathbb{F}_p,E)\to H_i(\mathbb{Z}/p^e\mathbb{Z},A,E)\to H_i(\mathbb{F}_p,A,E).$$
Then applying \ref{f-finit}, \ref{sirex} and \cite[Prop. IV.55]{A2} we get that $H_i(\mathbb{Z}/p^e\mathbb{Z},A,E)$ is finitely generated. In the same way, using the Jacobi-Zariski exact sequence in cohomology, we prove that the cohomology modules  $H^i(\mathbb{Z},A,E)$ and $H^i(\mathbb{Z}/p^e\mathbb{Z},A,E)$ are  finitely generated.

\begin{rem}\label{zrem}
 Looking  at Remark \ref{general} and Corollary \ref{z}, one can see that if $k$ is a perfect field (possibly of characteristic zero) and $A$ is a $k$-algebra which is a quotient of a regular ring, then $H_i(\mathbb{Z},A,E)$ is a finitely generated $A$-module for any $i\geq 1$ and any finitely generated $A$-module $E.$ 
\end{rem}

\begin{cor}\label{Zp}
 Let p be a prime number, A be a noetherian F-finite $\mathbb{F}_p$-algebra and E a finitely generated A-module. Then $H_i(\mathbb{Z}_{(p)},A,E)$ and $H^i(\mathbb{Z}_{(p)},A,E)$ are  finitely generated $A$-modules,
 for any $i\geq 0$.
\end{cor}

\par\noindent\textit{Proof:}  For $i=0$ we have $$H_0(\mathbb{Z}_{(p)},A,E)\simeq\Omega_{A/\mathbb{Z}_{(p)}}\otimes E\simeq \Omega_{A/\mathbb{Z}}\otimes E\simeq H_0(\mathbb{Z},A,E)$$ which is finitely generated by \ref{z}. For $i\geq 1,$ from the morphisms $\mathbb{Z}\to\mathbb{Z}_{(p)}\to A$ we have the exact sequence 
$$0=H_i(\mathbb{Z},\mathbb{Z}_{(p)},E)\to H_i(\mathbb{Z},A,E)\to H_i(\mathbb{Z}_{(p)},A,E)\to H_{i-1}(\mathbb{Z},\mathbb{Z}_{(p)},E)=0$$
and we apply again \ref{z}. The proof for the cohomology modules is similar.

\par We shall consider now the case of a $\mathbb{Z}_{(p)}$-algebra.
\begin{prop}\label{55}
 Let A be a noetherian F-finite $\mathbb{Z}_{(p)}$-algebra and E a finitely generated A-module. Then $H_i(\mathbb{Z}_{(p)},A,E)$ and $H^i(\mathbb{Z}_{(p)},A,E)$ are finitely 
 generated $A$-modules, for all $i\geq 1.$
\end{prop}

\par\noindent\textit{Proof:}  The case $pA=0,$ that is $A$ contains $\mathbb{F}_p,$ was considered above.
\par If $pA\neq 0,$ consider first the morphisms $\mathbb{Z}_{(p)}\to \mathbb{F}_p\to A/pA.$ 
We have the Jacobi-Zariski exact sequence 
$$H_i(\mathbb{Z}_{(p)},\mathbb{F}_p,E)\to H_i(\mathbb{Z}_{(p)},A/pA,E)\to H_i(\mathbb{F}_p,A/pA,E).$$
But $H_i(\mathbb{Z}_{(p)},\mathbb{F}_p,E)$ is finitely generated by \cite{A2}, IV.55 and $H_i(\mathbb{F}_p,A/pA,E)$ is finitely generated by \ref{f-finit}. By \ref{sirex} we obtain that $ H_i(\mathbb{Z}_{(p)},A/pA,E)$ is finitely generated. 
Consider now the morphisms $\mathbb{Z}_{(p)}\to A\to A/pA.$ 
We have the exact sequence 
$$H_{i+1}(A,A/pA,E)\to H_i(\mathbb{Z}_{(p)},A,E)\to H_i(\mathbb{Z}_{(p)},A/pA,E).$$
By the previous assertion $H_i(\mathbb{Z}_{(p)},A/pA,E)$ is finitely generated and by \cite{A2}, IV.55 $H_{i+1}(A,A/pA,E)$ is finitely generated. Now apply \ref{sirex} to get the assertion. 
\par For the cohomology modules, consider first the morphisms $\mathbb{Z}_{(p)}\to \mathbb{F}_p\to A/pA.$ 
We have the exact sequence 
$$H^i(\mathbb{F}_p,A/pA,E)\to H^i(\mathbb{Z}_{(p)},A/pA,E)\to H^i(\mathbb{Z}_{(p)},\mathbb{F}_p,E).$$
Then $H_i(\mathbb{Z}_{(p)},\mathbb{F}_p,E)$ is finitely generated by \cite{A2}, IV.55 and $H_i(\mathbb{F}_p,A/pA,E)$ is finitely generated by \ref{f-finit}. By \ref{sirex} we obtain that $ H_i(\mathbb{Z}_{(p)},A/pA,E)$ is finitely generated. 
Consider now the morphisms $\mathbb{Z}_{(p)}\to A\to A/pA.$ We have the associated Jacobi-Zariski exact sequence
$$H^i(\mathbb{Z}_{(p)},A/pA,E)\to H^i(\mathbb{Z}_{(p)},A,E)\to H^{i+1}(A,A/pA,E).$$
As in the proof of \ref{55} it follows that $H^i(\mathbb{Z}_{(p)},A,E)$ is finitely generated.

\begin{cor}\label{casez}
 Let A be a noetherian F-finite $\mathbb{Z}_{(p)}$-algebra and E a finitely generated A-module. Then $H_i(\mathbb{Z},A,E)$ and $H^i(\mathbb{Z},A,E)$ are finitely 
 generated $A$-modules, for all $i\geq 0.$
\end{cor}

\par\noindent\textit{Proof:}  It follows at once from \ref{55} and the Zariski-Jacobi exact sequence associated to 
$\mathbb{Z}\to\mathbb{Z}_{(p)}\to A$, taking account of the fact that $H_i(\mathbb{Z},\mathbb{Z}_{(p)},E)=H^i(\mathbb{Z},\mathbb{Z}_{(p)},E)=(0),\
\forall i\geq 0,$ cf. \cite[Prop. V.25]{A2}.

\bibliographystyle{amsplain}

\begin{thebibliography}{99}

\bibitem{A1} M. Andr\'e - \textit{M\'ethode simpliciale en alg\`ebre homologique et  alg\`ebre commutative} - Lecture Notes in Math. 32, Springer Verlag, Berlin Heidelberg New-York, 1967.

\bibitem{A2} M. Andr\'e - \textit{Homologie des alg\`ebres commutatives} - Springer Verlag, Berlin Heidelberg New-York, 1974.

\bibitem{A3} M. Andr\'e - \textit{Localisation de la lissit\'e formelle} - Manuscripta Math.,  13(1974), 297-307.

\bibitem{DM} B. I. Dundas, M. Morrow - \textit{Finite generation and continuity of topological Hochschild
and cyclic homology} - Ann. Scient. \'Ec. Norm. Sup.,  50(2017), 201-238.
\bibitem{Gab} O. Gabber - \textit{Notes on some t-structures} - in: Geometric Aspects of Dwork Theory, II, A. Adolphson, F. Baldassarri, P. Berthelot, N. Katz, F. Loeser (eds.) - Walter de
Gruyter GmbH \& Co. KG, Berlin, 2004, 711-734.
\bibitem{Iy} S. Iyengar - \textit{Andr\'e-Quillen homology of commutative algebras} - in: Interactions between Homotopy Theory and Algebra, L. L. Avramov, J. D. Christensen, W. G. Dwyer, M. A. Mandell, B. E. Shipley (eds.) - Contemporary Math., 436, 2007, 203-234.


\bibitem{LZ} A.  Langer, T. Zink - \textit{De Rham-Witt cohomology for a proper and smooth
morphism} - J. Inst. Math. Jussieu.,  3(2004), 231-314.

\bibitem{Qu} D. Quillen - \textit{On the(co-)homology of commutative rings} - in: Applications of categorical algebra, New York, 1968, Proc. Symp. Pure Math.17, Amer. Math. Soc., Providence, RI, 1970, 65–87.

\bibitem{Ro} J. Rotman - \textit{An Introduction to Homological Algebra, 2nd edition} - Springer, 2009.

\end{thebibliography}

\end{document}